# GLOBAL STABILIZATION OF FEEDFORWARD SYSTEMS UNDER PERTURBATIONS IN SAMPLING SCHEDULE

**Iasson Karafyllis[*] and Miroslav Krstic[**]**

[*]Dept. of Environmental Eng., Technical University of Crete,
73100, Chania, Greece, email: ikarafyl@enveng.tuc.gr

[**]Dept. of Mechanical and Aerospace Eng., University of California,
San Diego, La Jolla, CA 92093-0411, U.S.A., email: krstic@ucsd.edu

**Abstract**

For nonlinear systems that are known to be globally asymptotically stabilizable, control over networks introduces a major challenge because of the asynchrony in the transmission schedule. Maintaining global asymptotic stabilization in sampled-data implementations with zero-order hold and with perturbations in the sampling schedule is not achievable in general but we show in this paper that it is achievable for the class of feedforward systems. We develop sampled-data feedback stabilizers which are not approximations of continuous-time designs but are discontinuous feedback laws that are specifically developed for maintaining global asymptotic stabilizability under any sequence of sampling periods that is uniformly bounded by a certain "maximum allowable sampling period".

## 1. Introduction

Achieving stabilization by sampled-data feedback and ensuring robustness to perturbations in the sampling schedule are the central challenges in nonlinear control over networks, where asynchrony is ubiquitous. In this paper we achieve these goals for the class of uncertain feedforward systems, for which these goals are achievable due to the absence of exponential and finite escape time instabilities, despite the presence of nonlinearities of superlinear growth. We propose a saturation-based forwarding feedback, which we construct specifically for the sampled-data problem (namely, not as an approximation of a continuous design) and in such a way that it guarantees robustness of global asymptotic stability to perturbations in the sampling schedule.

Research on feedforward systems has played an important role in the development of nonlinear control theory, starting with the introduction of this class and the first feedback laws by Teel [54], followed by the key advances by Praly and Mazenc [31] and Jankovic, Sepulchre, and Kokotovic [15], and continuing with various extensions and generalizations by many authors [1-4,7,8,12,13,14,16,24,28-30,32-40,46-53,55-57]. More recently, feedforward systems with input delays and/or measurement delays have been studied [5,6,21,25].

In this work we focus on the problem of robust global stabilization of uncertain feedforward systems of the form

$$\begin{aligned}
\dot{x}_1 &= u \\
\dot{x}_2 &= x_1 + g_2(d, x_1, u) \\
&\vdots \\
\dot{x}_{n-1} &= x_{n-2} + g_{n-1}(d, x_1, \ldots, x_{n-2}, u) \\
\dot{x}_n &= x_{n-1} + g_n(d, x_1, \ldots, x_{n-1}, u) \\
x &= (x_1, \ldots, x_n)' \in \Re^n, u \in \Re, d \in D
\end{aligned} \quad (1.1)$$

where $D \subset \Re^l$ is a non-empty compact set and all mappings $g_i : D \times \Re^{i-1} \times \Re \to \Re$ ($i = 2, \ldots, n$) are locally Lipschitz for which there exists a smooth non-decreasing function $L \in C^0(\Re^+; \Re^+)$ such that

$$|g_i(d, \underline{x}_{i-1}, u)| \leq L(|(\underline{x}_{i-1}, u)|) |\underline{x}_{i-1}|^2 + L(|(\underline{x}_{i-1}, u)|) |\underline{x}_{i-1}| |u|, \ \underline{x}_{i-1} := (x_1, \ldots, x_{i-1})' \in \Re^{i-1}$$
$$\text{for all } (d, x, u) \in D \times \Re^n \times \Re \text{ and } i = 2, \ldots, n \quad (1.2)$$

More specifically, we solve the problem of robust global stabilization of (1.1) by means of bounded sampled-data feedback control applied with zero-order hold, i.e., with a controller of the form

$$\begin{aligned}
u(t) &= k(x(\tau_i)), \ t \in [\tau_i, \tau_{i+1}) \\
\tau_{i+1} &= \tau_i + r \exp(-w(\tau_i)), \ \tau_0 = 0 \\
w(t) &\in \Re^+
\end{aligned} \quad (1.3)$$

where $r > 0$ is a constant (the Maximum Allowable Sampling Period; MASP) and $k : \Re^n \to \Re$ is a bounded function with $k(0) = 0$. The input $w : \Re^+ \to \Re^+$ represents possible perturbations in the sampling schedule.

This problem is important for real-time implementations of control of feedforward systems, especially over networks, and to our knowledge, has not been addressed so far. The literature on sampled-data control provides control design methodologies that guarantee global stability for the following cases:

(i) Linear stabilizable systems, where $f(x,u) = Ax + Bu$, $A \in \Re^{n \times n}, B \in \Re^{n \times m}$,
(ii) Nonlinear systems of the form $\dot{x} = f(x) + g(x)u$, $x \in \Re^n, u \in \Re$, where the vector field $f : \Re^n \to \Re^n$ is globally Lipschitz and the vector field $g : \Re^n \to \Re^n$ is locally Lipschitz and bounded, which can be stabilized by a globally Lipschitz feedback law $u = k(x)$ (see [11]).
(iii) Nonlinear systems of the form $\dot{x}_i = f_i(x,u) + g_i(x,u)x_{i+1}$ for $i = 1, \ldots, n-1$ and $\dot{x}_n = f_n(x,u) + g_n(x,u)u$, where the drift terms $f_i(x,u)$ ($i = 1, \ldots, n$) satisfy the linear growth conditions $|f_i(x)| \leq L|x_1| + \ldots + L|x_i|$ ($i = 1, \ldots, n$) for certain constant $L \geq 0$ and there exist constants $b \geq a > 0$ such that $a \leq g_i(x,u) \leq b$ for all $i = 1, \ldots, n$, $x \in \Re^n, u \in \Re$ (see [19]).
(iv) Asymptotically controllable homogeneous systems with positive minimal power and zero degree (see [9]).
(v) Systems satisfying the reachability hypotheses of Theorem 3.1 in [20], or hypotheses (A1), (A2), (A3) in Section 4 of [18],
(vi) Nonlinear systems $\dot{x} = f(x,u)$, for which there exists a global diffeomorphism $\Theta : \Re^n \to \Re^n$ such that the change of coordinates $z = \Theta(x)$ transforms the system to one of the above cases.

However, nonlinear feedforward systems of the form (1.1) under hypothesis (1.2) rarely belong to one of the above classes (an exception is the class of linearizable feedforward systems; see [24,49-52]). On the other hand, there are well-established standard control design methodologies that guarantee stabilization of system (1.1) under sampled-data control with zero-order hold [10,22,23,26,27,41-45] but only in the practical and semiglobal sense. Therefore, the problem of robust global stabilization of (1.1) by means of bounded sampled-data feedback control applied with zero-order hold is open. Moreover, it was recently shown that the combination of a robust global sampled-data stabilizer with predictor schemes achieves global stabilization for systems with input and measurement delays [21]. Consequently, the solution of the problem of robust global stabilization of (1.1) by means of bounded sampled-data feedback control applied with zero-order hold automatically yields the solution of the same problem even in the presence of arbitrary measurement and input delays.

The key result of the present work is the "Sampled-Data Forwarding Lemma" (Lemma 3.1 below). The Sampled-Data Forwarding Lemma deals with a composite system that consists of two subsystems, the $x$-subsystem

$$\dot{x} = F(d,x,u)$$
$$x \in \Re^n, d \in D, u \in \Re$$
(1.4)

and the scalar $y$-subsystem

$$\dot{y} = G(d,x,u), y \in \Re \qquad (1.5)$$

Assuming that the $x$-subsystem can be robustly globally stabilized by means of sampled-data control applied with zero-order hold, then under appropriate conditions on the mappings $F,G$ we are in a position to construct a stabilizing feedback for the interconnected system (1.4), (1.5). The intuition behind the Sampled-Data Forwarding Lemma is as follows: we would like to bring the state $(x,y)$ of system (1.4), (1.5) to a neighborhood of the origin, where the linearization of (1.4), (1.5) prevails, and keep it there using a linear feedback strategy. In order to achieve this objective, we first apply the sampled-data feedback stabilizer for (1.4), which brings the $x$-component of the state $(x,y)$ close to zero. Once we have brought $x$ close to zero, we keep $x$ close to zero while simultaneously driving $y$ close to zero. Having brought $(x,y)$ to an appropriate neighborhood of zero, we apply linear feedback to drive $(x,y)$ to zero.

The Sampled-Data Forwarding Lemma provides an explicit formula for the robust feedback stabilizer and can be applied recursively for the robust global sampled-data stabilization of feedforward systems (Theorem 3.7 below). Moreover, if the assumed feedback stabilizer for (1.4) is bounded then the constructed feedback stabilizer for (1.4), (1.5) is bounded too. Robustness to perturbations in the sampling schedule is guaranteed by treating $x(\tau_i)$, where $\tau_i$ is a sampling time, as a perturbation of the current value of the state $x(t)$: by restricting the MASP $r > 0$, we are in a position to guarantee that $|x(t) - x(\tau_i)|$ is sufficiently small. The same methodology was introduced in the first author's papers [19,20], where robustness to perturbations in the sampling schedule and global stabilization was achieved for certain classes of nonlinear systems.

The structure of the paper is as follows: Section 2 provides the stability notions used in the paper and some technical results. Section 3 contains the "Sampled-Data Forwarding Lemma" (Lemma 3.1 below), which is applied recursively for the stabilization of (1.1). The main result (Theorem 3.7) guarantees the solvability of the problem of robust global stabilization of (1.1) by means of bounded sampled-data feedback control applied with zero-order hold. The formulae for the feedback stabilizers for feedforward systems contain parameters which can be tuned in order to guarantee good performance. A three-dimensional feedforward example is presented in Section 4 of the paper, which shows the importance of proper selection of the values of the parameters. Moreover, an additional example in Section 4 indicates that the Sampled-Data Forwarding Lemma is not restricted to feedforward systems. Section 5 contains the concluding remarks of the paper. Finally, the Appendix contains the proofs of all technical lemmas appearing in Section 3.

**Notation** Throughout this paper we adopt the following notation:

∗ Let $I \subseteq \Re_+ := [0,+\infty)$ be an interval. By $L^\infty(I;U)$ ($L^\infty_{loc}(I;U)$) we denote the space of measurable and (locally) essentially bounded functions $u(\cdot)$ defined on $I$ and taking values in $U \subseteq \Re^m$.

∗ By $C^0(A;\Omega)$, we denote the class of continuous functions on $A \subseteq \Re^n$, which take values in $\Omega \subseteq \Re^m$. By $C^k(A;\Omega)$, where $k \geq 1$, we denote the class of continuous functions on $A \subseteq \Re^n$, which have continuous derivatives of order $k \geq 1$ and take values in $\Omega \subseteq \Re^m$.

∗ For a vector $x \in \Re^n$ we denote by $x'$ its transpose and by $|x|$ its Euclidean norm. $A' \in \Re^{n \times m}$ denotes the transpose of the matrix $A \in \Re^{m \times n}$.

∗ By $sat: \Re \to [-1,1]$, we denote the continuous function $sat(x) = \dfrac{x}{\max(1,|x|)}$ for all $x \in \Re$. $B(x,\rho) \subseteq \Re^n$ denotes the closed ball in $\Re^n$ of radius $\rho \geq 0$ centered at $x \in \Re^n$, i.e., $B(x,\rho) := \{y \in \Re^n : |y - x| \leq \rho\}$.

∗ We say that an increasing continuous function $\gamma : \Re^+ \to \Re^+$ is of class $K_\infty$ if $\gamma(0) = 0$ and $\lim\limits_{s \to +\infty} \gamma(s) = +\infty$. By $KL$ we denote the set of all continuous functions $\sigma = \sigma(s,t) : \Re^+ \times \Re^+ \to \Re^+$ with the properties: (i) for each $t \geq 0$ the mapping $\sigma(\cdot,t)$ is non-decreasing; (ii) for each $s \geq 0$, the mapping $\sigma(s,\cdot)$ is non-increasing with $\lim\limits_{t \to +\infty} \sigma(s,t) = 0$.

## 2. Background Material and Preliminary Results

The stability notions used in the present work are applied to sampled-data systems of the form:

$$\begin{aligned} &\dot{x}(t) = f(d(t), x(t), x(\tau_i)) \quad , \quad t \in [\tau_i, \tau_{i+1}) \\ &\tau_0 = 0, \tau_{i+1} = \tau_i + r\exp(-w(\tau_i)), i = 0,1,... \\ &x(t) \in \Re^n, d(t) \in D, w(t) \in \Re^+ \end{aligned} \quad (2.1)$$

where $D \subset \Re^l$ is a non-empty set and $r > 0$ is a constant, under the following hypotheses:

**(H)** $f(d,x,x_0)$ is continuous with respect to $(d,x) \in D \times \Re^n$ and there exist a symmetric positive definite matrix $P \in \Re^{n \times n}$ and a function $a \in K_\infty$ such that the following inequalities hold

$$\sup\left\{\dfrac{(x-y)' P(f(d,x,x_0) - f(d,y,x_0))}{|x-y|^2} : x,y,x_0 \in B(0,s), x \neq y, d \in D\right\} < +\infty \, , \, \forall s > 0 \quad (2.2)$$

$$|f(d,x,x_0)| \leq a(|x| + |x_0|), \, \forall (u,d,x) \in U \times D \times \Re^n \quad (2.3)$$

Hypothesis (H) guarantees that $0 \in \Re^n$ is an equilibrium point for (1.1) and is automatically satisfied if $D \subset \Re^l$ is compact, $f(d,x,u)$ is locally Lipschitz with respect to $x \in \Re^n$ and $f(d,0,0) = 0$ for all $d \in D$. Moreover, by virtue of Proposition 2.5 in [17], hypothesis (H) guarantees that for every $(x_0, d, w) \in \Re^n \times L^\infty_{loc}(\Re^+; D) \times L^\infty_{loc}(\Re^+; \Re^+)$, system (2.1) admits a unique solution $x : [0, t_{\max}) \to \Re^n$ with $x(0) = x_0$, where $t_{\max} \in (0, +\infty]$ is the maximal existence time of the solution. Furthermore, if $t_{\max} < +\infty$ then $\limsup_{t \to t_{\max}^-} |x(t)| = +\infty$. The unique solution of (2.1) with $x(0) = x_0$ corresponding to inputs $(d,w) \in L^\infty_{loc}(\Re^+; D) \times L^\infty_{loc}(\Re^+; \Re^+)$ will be denoted by $x(t, x_0, d, w)$. The set of times $\{\tau_i\}_{i=0}^\infty$ is called the set of sampling times.

We next provide the definition of Robust Global Asymptotic Stability of (1.1).

**Definition 2.1:** *Consider system (2.1) under hypothesis (H). We say that $0 \in \Re^n$ is **Robustly Globally Asymptotically Stable (RGAS)** for system (2.1) if the following properties hold:*

**P1** (2.1) is **Robustly Lagrange Stable**, i.e., *for every $\varepsilon > 0$, it holds that*

$$\sup\{|x(t, x_0, d, w)| ; t \geq 0, |x_0| \leq \varepsilon, (d,w) \in L^\infty_{loc}(\Re^+; D \times \Re^+)\} < +\infty$$
**(Robust Lagrange Stability)**

**P2** (2.1) is **Robustly Lyapunov Stable**, i.e., *for every $\varepsilon > 0$ there exists $\delta := \delta(\varepsilon) > 0$ such that:*

$$|x_0| \leq \delta \Rightarrow |x(t, x_0, d, w)| \leq \varepsilon, \forall t \geq 0, \forall (d,w) \in L^\infty_{loc}(\Re^+; D) \times L^\infty_{loc}(\Re^+; \Re^+)$$
**(Robust Lyapunov Stability)**

**P3** (2.1) satisfies the **Robust Attractivity Property**, i.e. *for every $\varepsilon > 0$ and $R \geq 0$, there exists a $\tau := \tau(\varepsilon, R) \geq 0$, such that:*

$$|x_0| \leq R \Rightarrow |x(t, x_0, d, w)| \leq \varepsilon, \forall t \geq \tau, \forall (d,w) \in L^\infty_{loc}(\Re^+; D) \times L^\infty_{loc}(\Re^+; \Re^+)$$

**Remark 2.2:** Using Lemma 2.17 in [18] (with zero gain function) we can guarantee that $0 \in \Re^n$ is RGAS for system (2.1) if and only if there exists a function $\sigma \in KL$ such that the following estimate holds for all $(x_0, d, w) \in \Re^n \times L^\infty_{loc}(\Re^+; D) \times L^\infty_{loc}(\Re^+; \Re^+)$ and $t \geq 0$:

$$|x(t)| \leq \sigma(|x_0|, t) \qquad (2.4)$$

The reader should also notice that the sampling period is allowed to be time-varying. The factor $\exp(-w(\tau_i)) \leq 1$, with $w(t) \geq 0$ some non-negative function of time, is an uncertainty factor in the end-point of the sampling interval. Proving RGAS for (2.1) guarantees stability for all sampling schedules with $\tau_{i+1} - \tau_i \leq r$ (robustness to perturbations of the sampling schedule). Therefore, it is justified to call the constant $r > 0$ the Maximum Allowable Sampling Period (MASP).

We finish this section by providing a technical result which will be used in the following sections.

**Lemma 2.3:** *Let $b > a$ be constants and let $x : [a,b] \to \Re^n$ be absolutely continuous. Suppose that there exist constants $Q, G \geq 0$ such that*

$$|\dot{x}(t)| \leq Q|x(t)| + G|x(a)| \text{ for } t \in [a,b), \text{ a.e.} \tag{2.5}$$

*Suppose furthermore, that* $(G+Q)(b-a)\exp(Q(b-a)) < 1$. *Then the following inequality holds for all* $t \in [a,b]$:

$$|x(t) - x(a)| \leq \frac{(G+Q)(b-a)\exp(Q(b-a))}{1-(G+Q)(b-a)\exp(Q(b-a))}|x(t)| \tag{2.6}$$

**Proof:** Since $x:[a,b] \to \Re^n$ is absolutely continuous, it holds that $|x(t)-x(a)| \leq \int_a^t |\dot{x}(s)|ds$, for all $t \in [a,b]$. Inequality (2.5) implies $|x(t)-x(a)| \leq Q\int_a^t |x(s)|ds + G(b-a)|x(a)|$, for all $t \in [a,b]$ and consequently we obtain $|x(t)-x(a)| \leq Q\int_a^t |x(s)-x(a)|ds + (G+Q)(b-a)|x(a)|$, for all $t \in [a,b]$. Applying the Gronwall-Bellman lemma to the previous inequality, gives:

$$|x(t) - x(a)| \leq (G+Q)(b-a)\exp(Q(b-a))|x(a)|, \text{ for all } t \in [a,b]$$

The above inequality in conjunction with the triangle inequality implies that

$$|x(t) - x(a)| \leq (G+Q)(b-a)\exp(Q(b-a))|x(a)-x(t)| + (G+Q)(b-a)\exp(Q(b-a))|x(t)|, \forall t \in [a,b]$$

Since $(G+Q)(b-a)\exp(Q(b-a)) < 1$, the above inequality directly implies (2.6). The proof is complete. ◁

## 3. Main Results

All the results of the present work are proved by using Lemma 3.1, which is stated next. We call it the "Sampled-Data Forwarding Lemma" because it provides sufficient conditions for robust global stabilization by means of sampled-data control with positive sampling rate of a system with "added integration". The main result of the paper, Theorem 3.7, is established by recursively applying this lemma to system (1.1), and by constructively satisfying the lemma's key assumptions (inequalities (3.3)-(3.5) below) with the help of Lemma 3.6.

**Lemma 3.1 (The Sampled-Data Forwarding Lemma):** *Consider the following system*

$$\begin{aligned} \dot{x} &= Ax + bu + f(d,x,u) \\ \dot{y} &= x_n + g(d,x,u) \\ x &= (x_1,...,x_n)' \in \Re^n, y \in \Re, u \in \Re, d \in D \subset \Re^l \end{aligned} \tag{3.1}$$

*where* $b = (1,0,...,0)' \in \Re^n$, $D \subset \Re^l$ *is a nonempty compact set,* $f: D \times \Re^n \times \Re \to \Re^n$, $g: D \times \Re^n \times \Re \to \Re$ *are locally Lipschitz mappings with* $f(d,0,0) = 0$, $g(d,0,0) = 0$ *for all* $d \in D$, $A = \{a_{i,j} : i,j = 1,...,n\}$ *with* $a_{i,j} = 1$ *if* $j = i-1$, $i = 2,...,n$ *and* $a_{i,j} = 0$ *if otherwise. We assume that the following hypothesis holds:*

**(H)** *There exist a constant $r > 0$, a locally bounded function $k: \Re^n \to \Re$ with $k(0) = 0$ being continuous at $0 \in \Re^n$ such that $0 \in \Re^n$ is RGAS for the following sampled-data system:*

$$\begin{aligned}
\dot{x}(t) &= Ax(t) + bu(t) + f(d(t), x(t), u(t)) \\
u(t) &= k(x(\tau_i)), \, t \in [\tau_i, \tau_{i+1}) \\
\tau_{i+1} &= \tau_i + r\exp(-w(\tau_i)), \, \tau_0 = 0 \\
d(t) &\in D, w(t) \in \Re^+
\end{aligned} \quad (3.2)$$

*Let $P \in \Re^{n \times n}$ be a symmetric positive definite matrix, and $p \in \Re^n$ be a constant vector such that the matrix $P(A + bp') + (A' + pb')P$ is negative definite. Define $c = -(A' + pb')^{-1}(0,...,0,1)'$ and assume the existence of constants $M, R, K, \omega, \delta > 0$ such that*

$$\max\{x'P(Ax + f(d,x,u) + bu): (d,x) \in D \times \Re^n, x'Px = R^2, |u - p'x| \le K|c'b|\} < 0 \quad (3.3)$$

$$\max\{|g(d,x,u) + c'f(d,x,u)|: x'Px \le R^2, d \in D, |u - p'x| \le K|c'b|\} < K(c'b)^2 \quad (3.4)$$

$$z(Mg(d,x,u) + Mc'f(d,x,u) - K(c'b)\omega b'Px) \le (MK(c'b)^2\omega - \delta)|z|^2 - x'P((A + bp' + \delta I)x + f(d,x,u))$$
*for all $(x,z) \in \Re^n \times \Re$ with $x'Px \le R^2$, $\omega|z| \le 1$ and $u = p'x - Kc'b\omega z$* (3.5)

*Define $\tilde{k}: \Re^n \times \Re \to \Re$ by*

$$\tilde{k}(x,y) := \begin{cases} k(x) & \text{if } x'Px \ge R^2 \\ p'x - Kc'b\,\text{sat}(\omega(y + c'x)) & \text{if } x'Px < R^2 \end{cases} \quad (3.6)$$

*Then for sufficiently small $\tilde{r} > 0$, $0 \in \Re \times \Re^n$ is RGAS for the following sampled-data system*

$$\begin{aligned}
\dot{x}(t) &= Ax(t) + f(d(t), x(t), u(t)) + bu(t) \\
\dot{y}(t) &= x_n(t) + g(d(t), x(t), u(t)) \\
u(t) &= \tilde{k}(x(\tau_i), y(\tau_i)), \, t \in [\tau_i, \tau_{i+1}) \\
\tau_{i+1} &= \tau_i + \tilde{r}\exp(-w(\tau_i)), \, \tau_0 = 0 \\
d(t) &\in D, w(t) \in \Re^+
\end{aligned} \quad (3.7)$$

The proof of the Sampled-Data Forwarding Lemma is technical and is based on the following four technical results. Their proofs are provided in the Appendix.

**Lemma 3.2:** *Let $P \in \Re^{n \times n}$ be a symmetric positive definite matrix and $p \in \Re^n$ be a constant vector such that the matrix $P(A + bp') + (A' + pb')P$ is negative definite. Define $c = -(A' + pb')^{-1}(0,...,0,1)'$ and assume the existence of constants $M, R, K, \omega, \delta > 0$ such that (3.3), (3.4) and (3.5) hold. Consider the solution $(x(t), y(t)) \in \Re^n \times \Re$ of (3.7) under hypothesis (H), where $\tilde{k}: \Re^n \times \Re \to \Re$ is defined by (3.6) and $\tilde{r} > 0$, with arbitrary initial condition $(x(0), y(0)) \in \Re^n \times \Re$ satisfying $x'(0)Px(0) < R^2$ and corresponding to arbitrary $(d,w) \in L^\infty_{loc}(\Re^+; D \times \Re^+)$. If $\tilde{r} > 0$ is sufficiently small then $x'(t)Px(t) < R^2$ for all $t \in [0, \tau_1]$.*

Using induction and Lemma 3.2, we obtain the following result.

**Lemma 3.3:** Let $P \in \Re^{n \times n}$ be a symmetric positive definite matrix and $p \in \Re^n$ be a constant vector such that the matrix $P(A+bp')+(A'+pb')P$ is negative definite. Define $c = -(A'+pb')^{-1}(0,...0,1)'$ and assume the existence of constants $M, R, K, \omega, \delta > 0$ such that (3.3), (3.4) and (3.5) hold. Consider the solution $(x(t), y(t)) \in \Re^n \times \Re$ of (3.7) under hypothesis (H), where $\tilde{k} : \Re^n \times \Re \to \Re$ is defined by (3.6) and $\tilde{r} > 0$, with arbitrary initial condition $(x(0), y(0)) \in \Re^n \times \Re$ satisfying $x'(0)Px(0) < R^2$ and corresponding to arbitrary $(d, w) \in L^\infty_{loc}(\Re^+; D \times \Re^+)$. If $\tilde{r} > 0$ is sufficiently small then the solution $(x(t), y(t)) \in \Re^n \times \Re$ of (3.7) exists for all $t \geq 0$ and satisfies $x'(t)Px(t) < R^2$ for all $t \geq 0$.

The following lemma uses the result of Lemma 3.3 and shows attractivity for a certain region in the state space.

**Lemma 3.4:** Let $P \in \Re^{n \times n}$ be a symmetric positive definite matrix and $p \in \Re^n$ be a constant vector such that the matrix $P(A+bp')+(A'+pb')P$ is negative definite. Define $c = -(A'+pb')^{-1}(0,...0,1)'$ and assume the existence of constants $M, R, K, \omega, \delta > 0$ such that (3.3), (3.4) and (3.5) hold. Consider the solution $(x(t), y(t)) \in \Re^n \times \Re$ of (3.7) under hypothesis (H), where $\tilde{k} : \Re^n \times \Re \to \Re$ is defined by (3.6) and $\tilde{r} > 0$, with arbitrary initial condition $(x(0), y(0)) \in \Re^n \times \Re$ satisfying $x'(0)Px(0) < R^2$ and corresponding to arbitrary $(d, w) \in L^\infty_{loc}(\Re^+; D \times \Re^+)$. If $\tilde{r} > 0$ is sufficiently small then there exists $T \in C^0(\Re; \Re^+)$ such that

$$|z(t)| \leq \max\{|z(0)|, \omega^{-1}\}, \quad \forall t \geq 0 \tag{3.8}$$

$$|z(t)| \leq \omega^{-1}, \quad \forall t \geq T(z(0)) \tag{3.9}$$

where $z(t) = y(t) + c'x(t)$.

**Lemma 3.5:** Let $P \in \Re^{n \times n}$ be a symmetric positive definite matrix and $p \in \Re^n$ be a constant vector such that the matrix $P(A+bp')+(A'+pb')P$ is negative definite. Define $c = -(A'+pb')^{-1}(0,...0,1)'$ and assume the existence of constants $M, R, K, \omega, \delta > 0$ such that (3.3), (3.4) and (3.5) hold. Consider the solution $(x(t), y(t)) \in \Re^n \times \Re$ of (3.7) under hypothesis (H), where $\tilde{k} : \Re^n \times \Re \to \Re$ is defined by (3.6) and $\tilde{r} > 0$, with arbitrary initial condition $(x(0), y(0)) \in \Re^n \times \Re$ satisfying $x'(0)Px(0) < R^2$, $|y(0)+c'x(0)| \leq \omega^{-1}$ and corresponding to arbitrary $(d, w) \in L^\infty_{loc}(\Re^+; D \times \Re^+)$.

If $\tilde{r} > 0$ is sufficiently small then there exists $\mu > 0$ such that the following differential inequality holds for $t \geq 0$, a.e.:

$$\dot{V}(t) \leq -\mu V(t) \tag{3.10}$$

where $z(t) = y(t) + c'x(t)$ and $V(t) = \frac{M}{2}z^2(t) + \frac{1}{2}x'(t)Px(t)$.

**Remark:** The reader should notice that by virtue of Lemma 3.3 and Lemma 3.4 the set $S = \{(x, y) \in \Re^n \times \Re : x'Px < R^2, |y+c'x| \leq \omega^{-1}\}$ is positively invariant for system (3.7). The differential inequality (3.10) guarantees that $V(t) \leq \exp(-\mu t)V(0)$ for all $t \geq 0$, provided that $(x(0), y(0)) \in S$.

Since $V(x,y) = \frac{M}{2}|y+c'x|^2(t) + \frac{1}{2}x'Px$ is a positive definite quadratic function, the previous inequality shows that local exponential stability is guaranteed for system (3.7) in the region $S \subseteq \Re^n \times \Re$. Notice that the size of the region $S \subseteq \Re^n \times \Re$ is determined by the constants $R, \omega$.

We are now in a position to prove the Sampled-Data Forwarding Lemma.

**Proof of Lemma 3.1:** We will restrict $\tilde{r} > 0$, so that

$$\tilde{r} \leq r \tag{3.11}$$

Notice that Lemma 3.3 and definition (3.6) imply that the set $\{x \in \Re^n : x'Px < R^2\}$ is positively invariant. Robust Lyapunov stability for system (3.7) is a direct consequence of the differential inequality (3.10). We will show next robust Lagrange stability and robust attractivity for system (3.7).

Consider the solution $(x(t), y(t)) \in \Re^n \times \Re$ of (3.7) under hypotheses (H1-2), with initial condition $(x(0), y(0)) \in \Re^n \times \Re$ and corresponding to arbitrary $(d,w) \in L^\infty_{loc}(\Re^+; D \times \Re^+)$. By virtue of hypothesis (H), inequality (3.11) and definition (3.6), there exists $\sigma \in KL$ such that

$$|x(t)| \leq \sigma(|x(0)|, t) \tag{3.12}$$

for all times $t \geq 0$ with $x'(t)Px(t) \geq R^2$. Inequality (3.12) in conjunction with Lemma 3.3 implies that there exists a constant $C > 0$ such that the following inequality holds:

$$|x(t)| \leq \max(\sigma(|x(0)|, 0), C), \quad \forall t \geq 0 \tag{3.13}$$

and that there exists $\tilde{T} \in C^0(\Re^n; \Re^+)$ such that

$$x'(t)Px(t) < R^2, \text{ for all } t \geq \tilde{T}(x(0)) \tag{3.14}$$

Notice that hypothesis (H) implies the existence of $\rho \in K_\infty$ such that $|k(x)| \leq \rho(|x|)$, for all $x \in \Re^n$. Therefore, using (3.13), (3.14) we can conclude that there exists $\gamma \in K_\infty$ such that:

$$|y(t) + c'x(t)| \leq \gamma(|x(0)| + |y(0)|) \tag{3.15}$$

for all times $t \geq 0$ with $x'(t)Px(t) \geq R^2$. Therefore, inequality (3.15) in conjunction with Lemma 3.4 implies that:

$$|y(t) + c'x(t)| \leq \max(\gamma(|x(0)| + |y(0)|), |y(0) + c'x(0)|, \omega^{-1}), \quad \forall t \geq 0 \tag{3.16}$$

Estimates (3.13), (3.16) prove robust Lagrange stability. Finally, inequality (3.14) in conjunction with Lemma 3.4 and Lemma 3.5 imply that robust attractivity holds as well. The proof is complete.   ◁

The following result shows that the assumptions of the Sampled-Data Forwarding Lemma can be automatically satisfied for a certain class of nonlinearities.

**Lemma 3.6:** *Suppose that there exists a non-decreasing function $L \in C^0(\Re^+; \Re^+)$ such that the following inequality holds for the mappings $f : D \times \Re^n \times \Re \to \Re^n$, $g : D \times \Re^n \times \Re \to \Re$ :*

$$|f(d,x,u)| + |g(d,x,u)| \leq L(|(x,u)|)|x|^2 + L(|(x,u)|)|x||u|, \text{ for all } (d,x,u) \in D \times \Re^n \times \Re \quad (3.17)$$

*Then there exist constants $R^*, C > 0$ with $C \leq 1$ such that for every $\omega > 0$, $R \in (0, R^*)$ there exist constants $M, \delta > 0$ such that (3.3), (3.4), (3.5) hold with $K = CR$.*

**Theorem 3.7:** *Consider system (1.1) where all mappings $g_i : D \times \Re^{i-1} \times \Re \to \Re$ are locally Lipschitz and assume that there exists a smooth non-decreasing function $L \in C^0(\Re^+; \Re^+)$ such that (1.2) holds. Then there exist a bounded $k : \Re^n \to \Re$ with $k(0) = 0$ being continuous at $0 \in \Re^n$ and a constant $r > 0$, such that $0 \in \Re^n$ is RGAS for the closed-loop sampled-data system (1.1) with*

$$\begin{aligned} u(t) &= k(x(\tau_i)), \, t \in [\tau_i, \tau_{i+1}) \\ \tau_{i+1} &= \tau_i + r \exp(-w(\tau_i)), \, \tau_0 = 0 \\ w(t) &\in \Re^+ \end{aligned} \quad (3.18)$$

*Define $Q_i \in \Re^{i \times n}$ with $Q_i x = (x_1, \ldots x_i)'$ for $i = 1, \ldots, n$, $b_1 = [1] \in \Re$, $b_i = (1, 0, \ldots, 0)' \in \Re^i$ for $i = 2, \ldots, n$, $A_1 = [0] \in \Re^{1 \times 1}$, $A_i = \{a_{k,j} : k, j = 1, \ldots, i\} \in \Re^{i \times i}$ for $i = 2, \ldots, n$ with $a_{k,j} = 1$ if $j = k-1$, $k = 2, \ldots, i$ and $a_{k,j} = 0$ if otherwise. Let arbitrary constants $K_0 > 0$, $\omega_i > 0$ ($i = 0, \ldots, n-1$), arbitrary matrices $P_i \in \Re^{i \times i}$ ($i = 1, \ldots, n-1$) being symmetric and positive definite and arbitrary vectors $p_i \in \Re^i$ ($i = 1, \ldots, n-1$) such that the matrices $P_i (A_i + b_i p_i') + (A_i' + p_i b_i') P_i$ ($i = 1, \ldots, n-1$) are negative definite. Define $c_i = -(A_i' + p_i b_i')^{-1}(0, \ldots 0, 1)' \in \Re^i$ for $i = 1, \ldots, n-1$. Then there exist constants $r > 0$, $K_i > 0$, $R_i > 0$ ($i = 1, \ldots, n-1$) such that $0 \in \Re^n$ is RGAS for the closed-loop sampled-data system (1.1) with (3.18), where $k : \Re^n \to \Re$ is defined by*

$$k(x) := p_i' Q_i x - K_i c_i' b_i \operatorname{sat}(\omega_i (x_{i+1} + c_i' Q_i x)) \quad (3.19)$$

*where $i = i(x) \in \{1, \ldots, n-1\}$ is the largest integer such that*

$$x' Q_i' P_i Q_i x < R_i^2 \quad (3.20)$$

*and*

$$k(x) := -K_0 \operatorname{sat}(\omega_0 x_1), \text{ if } \min_{l=1,\ldots,n-1} \left( x' Q_l' P_l Q_l x - R_l^2 \right) \geq 0 \quad (3.21)$$

**Proof:** Repeated application of the Sampled-Data Forwarding Lemma and Lemma 3.6. Notice that the subsystem $\dot{x}_1 = u$ can be stabilized by the bounded sampled-data feedback

$$\begin{aligned} u(t) &= -K_0 \operatorname{sat}(x_1(\tau_i)), \, t \in [\tau_i, \tau_{i+1}) \\ \tau_{i+1} &= \tau_i + K_0^{-1} \exp(-w(\tau_i)), \, \tau_0 = 0 \\ w(t) &\in \Re^+ \end{aligned}$$

where $K_0 > 0$ is an arbitrary positive constant. The Sampled-Data Forwarding Lemma is applied for $j = 1,...,n-1$ with $x \in \Re^n$ replaced by $(x_1,...,x_j)' \in \Re^j$, $y \in \Re$ replaced by $x_{j+1} \in \Re$, $A \in \Re^{n \times n}$ replaced by $A_j \in \Re^{j \times j}$, $b \in \Re^n$ replaced by $b_j \in \Re^j$, $g(d,x,u) \in \Re$ replaced by $g_{j+1}(d,x_1,...,x_j,u) \in \Re$, $f(d,x,u) \in \Re^n$ replaced by $f(d,x,u) = (0, g_2(d,x_1,u),...,g_j(d,x_1,...,x_{j-1},u))' \in \Re^j$ for $j \geq 2$ and $f(d,x,u) = 0 \in \Re$ for $j = 1$, $P \in \Re^{n \times n}$, $p \in \Re^n$, $c = -(A' + pb')^{-1}(0,...0,1)'$ replaced by $P_j \in \Re^{j \times j}$, $p_j \in \Re^j$, $c_j = -(A'_j + p_j b'_j)^{-1}(0,...0,1)' \in \Re^j$, respectively and $k : \Re^n \to \Re$ replaced by $k_j : \Re^j \to \Re$ which is defined by the following equalities:

- for $j \geq 2$

$$k_j(x_1,...,x_j) := p'_i Q_i x - K_i c'_i b_i \, \text{sat}(\omega_i (x_{i+1} + c'_i Q_i x)) \tag{3.22}$$

where $i = i(x_1,...,x_{j-1}) \in \{1,...,j-1\}$ is the largest integer such that (3.20) holds and

$$k_j(x_1,...,x_j) := -K_0 \, \text{sat}(\omega_0 x_1), \text{ if } \min_{l=1,...,j-1}\left(x' Q'_l P_l Q_l x - R_l^2\right) \geq 0 \tag{3.23}$$

- for $j = 1$

$$k_1(x_1) := -K_0 \, \text{sat}(\omega_0 x_1) \tag{3.24}$$

By virtue of (1.2), it follows that (3.17) holds with $x \in \Re^n$ replaced by $(x_1,...,x_j)' \in \Re^j$, $g(d,x,u) \in \Re$ replaced by $g_{j+1}(d,x_1,...,x_j,u) \in \Re$, $f(d,x,u) \in \Re^n$ replaced by $f(d,x,u) = (0, g_2(d,x_1,u),...,g_j(d,x_1,...,x_{j-1},u))' \in \Re^j$ for $j \geq 2$ and $f(d,x,u) = 0 \in \Re$ for $j = 1$ and $L \in C^0(\Re^+; \Re^+)$ replaced by $L_j \in C^0(\Re^+; \Re^+)$, where $L_j(s) = jL(s)$ and $L$ is the function involved in (1.2). Therefore, by virtue of Lemma 3.6, it follows that there exist constants $R^*_j, C_j > 0$ with $C_j \leq 1$ such that for every $\omega_j > 0$, $R_j \in (0, R^*_j)$ there exist constants $M_j, \delta_j > 0$ such that (3.3), (3.4), (3.5) hold with $K = K_j = C_j R_j$, $R = R_j$, $\omega = \omega_j$, $M = M_j$ and $\delta = \delta_j$.

The proof is complete. ◁

**Remark 3.8:** Notice that the proof of Theorem 3.7 guarantees that for every $G > 0$, the sampled-data feedback stabilizer $k : \Re^n \to \Re$ can be selected in such a way that $|k(x)| \leq G$ for all $x \in \Re^n$. To see this, first select arbitrary constants $\omega_i > 0$ ($i = 0,...,n-1$), arbitrary matrices $P_i \in \Re^{i \times i}$ ($i = 1,...,n-1$) being symmetric and positive definite and arbitrary vectors $p_i \in \Re^i$ ($i = 1,...,n-1$) such that the matrices $P_i(A_i + b_i p'_i) + (A'_i + p_i b'_i)P_i$ ($i = 1,...,n-1$) are negative definite. The selection of $K_i > 0$, $R_i > 0$ ($i = 1,...,n-1$) made in the proof of Theorem 3.7 guarantees that the constants $R_i \in (0, R^*_i)$ ($i = 1,...,n-1$) can be selected in an arbitrary way, where $R^*_i > 0$ ($i = 1,...,n-1$) are appropriate constants. Moreover, the inequalities $K_i \leq R_i$ hold for $i = 1,...,n-1$. It follows from (3.19), (3.20), (3.21) that

$$|k(x)| \leq \max\left\{K_0, \max_{i=1,...,n-1} R_i \left(|p_i|a_i^{-1} + |c'_i b_i|\right)\right\}, \quad \forall x \in \Re^n \tag{3.25}$$

where $a_i > 0$ ($i = 1,...,n-1$) are constants satisfying $x'Q_i'P_iQ_ix \geq a_i^2|Q_ix|^2$ for all $x \in \Re^n$. It follows from (3.25) that if $K_0 \leq G$ and $R_i \leq \dfrac{G}{|p_i|a_i^{-1} + |c_i'b_i|}$ for $i = 1,...,n-1$ then $|k(x)| \leq G$ for all $x \in \Re^n$. Notice that we can always select $K_0 \leq G$ and $R_i \leq \dfrac{G}{|p_i|a_i^{-1} + |c_i'b_i|}$ for $i = 1,...,n-1$ ($K_0 > 0$ and $R_i \in (0, R_i^*)$ are free parameters).

## 4. Illustrative Examples

In this section we present two examples that illustrate the results of the previous section. The first example shows the application of Theorem 3.7 to a feedforward system.

**Example 4.1:** We consider the 3-dimensional feedforward system

$$\begin{aligned} \dot{x}_1 &= u \\ \dot{x}_2 &= x_1 + x_1 u \\ \dot{x}_3 &= x_2 + x_1^2 \\ x &= (x_1, x_2, x_3)' \in \Re^3, u \in \Re \end{aligned} \quad (4.1)$$

The solution map of system (4.1) can be explicitly found: the resulting discrete-time system that corresponds to a constant sampling period $r > 0$ and input $u \in \Re$ applied with zero order hold is given by the following equations:

$$\begin{aligned} x_1^+ &= x_1 + ur \\ x_2^+ &= x_2 + (x_1 + ux_1)r + (u + u^2)\dfrac{r^2}{2} \\ x_3^+ &= x_3 + (x_2 + x_1^2)r + (x_1 + 3ux_1)\dfrac{r^2}{2} + (u + 3u^2)\dfrac{r^3}{6} \end{aligned} \quad (4.2)$$

However, as already noted in the Introduction, system (4.1) is not included in one of the classes of systems noted in the Introduction for which there exists a feedback design methodology that results in the design of a globally stabilizing sampled-data feedback (notice that (4.1) is not linearizable). Other approaches for sampled-data systems can be also applied (see [10,22,23,26,27,41-45]) but the result is semiglobal and practical sampled-data stabilization of system (4.1).

Here we apply the step-by-step feedback design methodology described in Theorem 3.7. The feedback law will be given by (3.19), (3.20), (3.21). For simplicity, we select $\omega_0 = \omega_1 = \omega_2 = 1$ and $K_0 = 1$. We also select:

$$P_1 = [1], \quad p_1 = [-1], \quad P_2 = \begin{bmatrix} 1 & 1 \\ 1 & 2 \end{bmatrix}, \quad p_2 = \begin{bmatrix} -2 \\ -2 \end{bmatrix}$$

Using the formula $c_i = -(A_i' + p_ib_i')^{-1}(0,...0,1)' \in \Re^i$ for $i = 1,2$, we obtain $c_1 = [1]$ and $c_2 = \begin{bmatrix} 1/2 \\ 1 \end{bmatrix}$. The only constants that remain to be determined are $R_1, R_2, K_1, K_2$.

In order to determine $R_1, K_1 > 0$, we use the Sampled-Data Forwarding Lemma. We apply the Sampled-Data Forwarding Lemma with $n=1$, $A=[0]$, $b=[1]$, $P=[1]$, $p=[-1]$, $c=[1]$, $f(d,x,u) \equiv 0$ and $g(d,x,u) = x_1 u$. Conditions (3.3), (3.4), (3.5) are satisfied with $M = \frac{K}{R+K}$, $\omega = 1$ for $\delta > 0$ sufficiently small, provided that

$$\frac{R^2}{1-R} < K < R \quad \text{and} \quad R + K < 1 \tag{4.3}$$

Inequalities (4.3) hold with $R = \frac{3}{8}$ and $K = \frac{1}{4}$. Therefore, we select $R_1 = \frac{3}{8}$ and $K_1 = \frac{1}{4}$.

In order to determine $R_2, K_2 > 0$, we use again the Sampled-Data Forwarding Lemma. We apply the Sampled-Data Forwarding Lemma with $n=2$, $A = \begin{bmatrix} 0 & 0 \\ 1 & 0 \end{bmatrix}$, $b = \begin{bmatrix} 1 \\ 0 \end{bmatrix}$, $P = \begin{bmatrix} 1 & 1 \\ 1 & 2 \end{bmatrix}$, $p = \begin{bmatrix} -2 \\ -2 \end{bmatrix}$, $c = \begin{bmatrix} 1/2 \\ 1 \end{bmatrix}$, $f(d,x,u) = \begin{bmatrix} 0 \\ x_1 u \end{bmatrix}$ and $g(d,x,u) = x_1^2$. After some tedious calculations, we conclude that conditions (3.3), (3.4), (3.5) are satisfied with $M = K \frac{2+(3+2\sqrt{2})R}{4R}$, $\omega = 1$ for $\delta > 0$ sufficiently small, provided that

$$\frac{4R^2}{1-2\sqrt{2}R} < K < 2R \frac{1-2(2+\sqrt{2})R}{R+1} \quad \text{and} \quad (4+2\sqrt{2})R + (3-2\sqrt{2})R^2 < 1 \tag{4.4}$$

Inequalities (4.4) hold with $R = K = \frac{1}{20}$. Therefore, we select $R_2 = K_2 = \frac{1}{20}$. We conclude that the sampled-data feedback law (3.18) defined by

$$k(x) := \begin{cases} -\operatorname{sat}(x_1) & \text{if } 8|x_1| \geq 3 \text{ and } 20\sqrt{x_2^2 + (x_1+x_2)^2} \geq 1 \\ -x_1 - \frac{1}{4}\operatorname{sat}(x_2 + x_1) & \text{if } 8|x_1| < 3 \text{ and } 20\sqrt{x_2^2 + (x_1+x_2)^2} \geq 1 \\ -2(x_1+x_2) - \frac{1}{40}\operatorname{sat}\left(x_3 + x_2 + \frac{1}{2}x_1\right) & \text{if } 20\sqrt{x_2^2 + (x_1+x_2)^2} < 1 \end{cases} \tag{4.5}$$

achieves global stabilization of system (4.1) provided that the MASP $r > 0$ is sufficiently small. Indeed, simulations show that global stabilization of system (4.1) is achieved with $r = 0.01$. However, under these conservative choices of design parameters, which satisfy the sufficient conditions of Theorem 3.7, the closed-loop system shows different dynamic behavior in different time scales. The state variables $x_1, x_2$ converge very fast, while the state variable $x_3$ exhibits slow convergence, which lasts about 900 time units.

Therefore, it is crucial to determine tight bounds for the range of values for $R_2, K_2 > 0$ which guarantee global asymptotic stability. Numerical experiments show that higher values than 0.05 for $R_2, K_2 > 0$ can guarantee global asymptotic stability for $r = 0.2$. Figure 1 shows the evolution of the state variables for the closed-loop system (4.1) with (3.18), where $k$ is defined by

$$k(x) := \begin{cases} -\text{sat}(x_1) & \text{if} \quad 8|x_1| \geq 3 \quad \text{and} \quad \sqrt{x_2^2 + (x_1 + x_2)^2} \geq 1 \\ -x_1 - \dfrac{1}{4}\text{sat}(x_2 + x_1) & \text{if} \quad 8|x_1| < 3 \quad \text{and} \quad \sqrt{x_2^2 + (x_1 + x_2)^2} \geq 1 \\ -2(x_1 + x_2) - \dfrac{1}{2}\text{sat}\left(x_3 + x_2 + \dfrac{1}{2}x_1\right) & \text{if} \quad \sqrt{x_2^2 + (x_1 + x_2)^2} < 1 \end{cases} \quad (4.6)$$

with $r = 0.2$, $w(t) = \ln\left(\dfrac{2}{1 + |\sin(t)|}\right)$ and initial condition $x_1(0) = x_2(0) = x_3(0) = 1$. It is clear that the selection $R_2 = K_2 = 1$ guarantees good performance even when perturbations of the sampling schedule are present. Figure 2 shows the corresponding input behavior and Figure 3 focuses on the evolution of the input for $t \in [4,8]$.

Having addressed the robust global stabilization problem for (4.1) under sampled-data control applied with zero order hold, we are in a position to address the stabilization problem for the system

$$\dot{x}_1(t) = u(t - \tau), \quad \dot{x}_2(t) = x_1(t) + x_1(t)u(t - \tau), \quad \dot{x}_3(t) = x_2(t) + x_1^2(t) \qquad (4.7)$$
$$x = (x_1, x_2, x_3)' \in \Re^3, \quad u \in \Re$$

where $\tau > 0$ is the input delay, with sampled and delayed measurement

$$y(t) = x(\tau_i - T) \qquad (4.8)$$

where $\tau_i = ir, i \in Z^+$ are the sampling times and $r > 0$, $T > 0$ are the sampling period and the measurement delay, respectively. Following the methodology described in [], first we select a sampling period for which there exists $l \in Z^+$ such that $\tau = lr$. This is possible since we have shown robustness with respect to the sampling schedule, i.e., we have shown robust global asymptotic stability for the closed-loop system (4.1) with (3.18), where $k$ is defined by (4.6) with $r = 0.2$. Consequently, we may choose any integer $l \in Z^+$ with $l \geq 5\tau$ and set $r = \dfrac{\tau}{l} \leq 0.2$. The predictor mapping for (4.1) can be explicitly expressed by the equations:

$$\Phi(x,u) := \left[ x_1 + \int_{-T-\tau}^{0} u(s)ds, \quad x_2 + (\tau + T)x_1 + x_1\int_{-T-\tau}^{0} u(s)ds + \int_{-T-\tau}^{0}(1 + u(s))\int_{-T-\tau}^{s} u(q)dq\,ds, \quad \phi(x,u) \right]' \qquad (4.9)$$

where

$$\phi(x,u) = x_3 + (\tau + T)x_2 + (\tau + T)x_1^2 + \dfrac{1}{2}(\tau + T)^2 x_1 + 3x_1 \int_{-T-\tau}^{0}\int_{-T-\tau}^{s} u(q)dq\,ds$$
$$+ \int_{-T-\tau}^{0}\int_{-T-\tau}^{s}(1 + u(w))\int_{-T-\tau}^{w} u(q)dq\,dw\,ds + \int_{-T-\tau}^{0}\left(\int_{-T-\tau}^{s} u(q)dq\right)^2 ds \qquad (4.10)$$

Using Theorem 3.2 in [21], we can guarantee that the closed-loop system (4.7) with

$$u(t) = k(X_1(\tau_i), X_2(\tau_i), X_3(\tau_i)), \quad t \in [\tau_i, \tau_{i+1}), i \in Z^+ \qquad (4.11)$$
$$\tau_{i+1} = \tau_i + r, \tau_0 = 0$$

where $k$ is defined by (4.6) and

$$X_1(\tau_i) = x_1(\tau_i - T) + \int_{\tau_i - T - \tau}^{\tau_i} u(s)ds$$

$$X_2(\tau_i) = x_2(\tau_i - T) + (\tau + T)x_1(\tau_i - T) + x_1(\tau_i - T) \int_{\tau_i - T - \tau}^{\tau_i} u(s)ds + \int_{\tau_i - T - \tau}^{\tau_i} (1 + u(s)) \int_{\tau_i - T - \tau}^{s} u(q)dq\, ds$$

(4.12)

$$X_3(\tau_i) = x_3(\tau_i - T) + (\tau + T)\left(x_2(\tau_i - T) + x_1^2(\tau_i - T)\right) + \frac{1}{2}(\tau + T)^2 x_1(\tau_i - T)$$

$$+ 3x_1(\tau_i - T) \int_{\tau_i - T - \tau}^{\tau_i} \int_{\tau_i - T - \tau}^{s} u(q)dq\, ds + \int_{\tau_i - T - \tau}^{\tau_i} \int_{\tau_i - T - \tau}^{s} (1 + u(w)) \int_{\tau_i - T - \tau}^{w} u(q)dq\, dw\, ds + \int_{\tau_i - T - \tau}^{\tau_i} \left(\int_{\tau_i - T - \tau}^{s} u(q)dq\right)^2 ds$$

(4.13)

is Uniformly Globally Asymptotically Stable, in the sense that there exists a function $\sigma \in KL$ such that for every $(x_0, u_0) \in C^0([-T,0]; \Re^3) \times L^\infty([-T-\tau,0); \Re)$, the solution $(x(t), u(t)) \in \Re^3 \times \Re$ of system (4.7), (4.11) with initial condition $u(\theta) = u_0(\theta)$ for $\theta \in [-\tau - T, 0)$, $x(\theta) = x_0(\theta)$ for $\theta \in [-T, 0]$ satisfies the following inequality for all $t \geq 0$:

$$\max_{t-T \leq s \leq t} |x(s)| + \sup_{t-T-\tau \leq s < t} |u(s)| \leq \sigma\left(\max_{-T \leq s \leq 0} |x(s)| + \sup_{-T-\tau \leq s < 0} |u(s)|, t\right)$$

(4.14)

It should be emphasized that the value of the integrals involved in (4.12), (4.13) can be computed with precision when $\tau_i \geq T + \tau$, since the input $u(t)$ for $t \geq 0$ is a piecewise constant function.  ◁

The second example shows that the Sampled-Data Forwarding Lemma (Lemma 3.1) can be also applied to some nonlinear systems outside of the class of feedforward systems.

**Example 4.2:** Consider the nonlinear system

$$\begin{aligned}\dot{x} &= Ax + bu + f(d, x) \\ \dot{y} &= x_n + g(d, x) \\ x &\in \Re^n, y \in \Re, d \in D, u \in \Re\end{aligned}$$

(4.15)

where $b = (1, 0, ..., 0)' \in \Re^n$, $A = \{a_{i,j} : i, j = 1, ..., n\}$ with $a_{i,j} = 1$ if $j = i - 1$, $i = 2, ..., n$ and $a_{i,j} = 0$ if otherwise, $D \subset \Re^l$ is a non-empty compact set, $f : D \times \Re^n \to \Re^n$, $g : D \times \Re^n \to \Re$ are locally Lipschitz mappings that satisfy the following inequalities:

$$\max\{|f(d, x)|, |g(d, x)|\} \leq L_1 |x|, \text{ for all } (d, x) \in D \times B(0, \rho)$$

(4.16)

$$|f(d, x)| \leq L_2 |x|, \text{ for all } (d, x) \in D \times \Re^n$$

(4.17)

for certain constants $L_2 \geq L_1 > 0$ and $\rho > 0$. At this point we should note the crucial difference between (4.16), (4.17) and (3.17). While in (3.17) the nonlinearities $f$ and $g$ are restricted to be locally quadratic, in (4.16), (4.17) the nonlinearities are allowed to have linear growth at the origin.

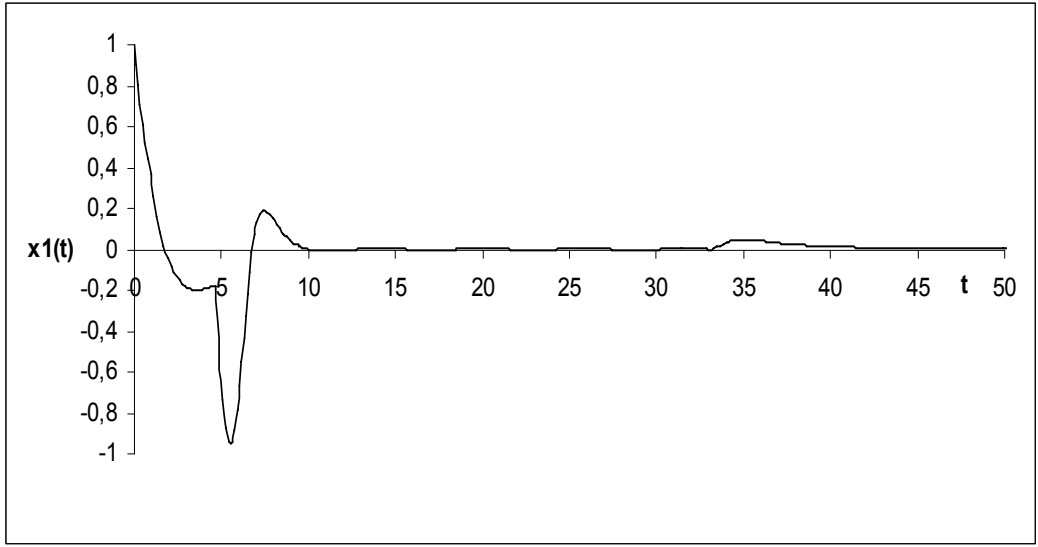

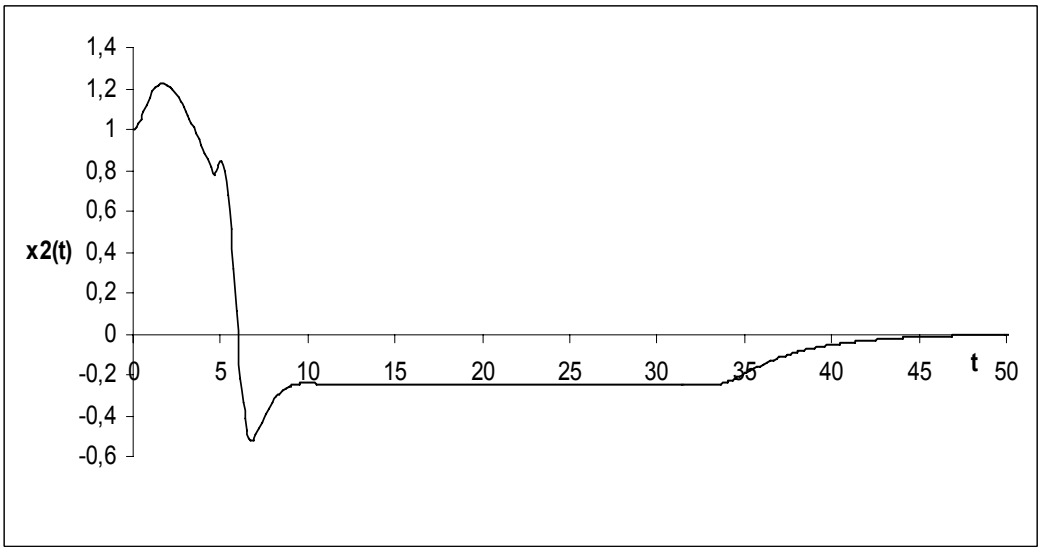

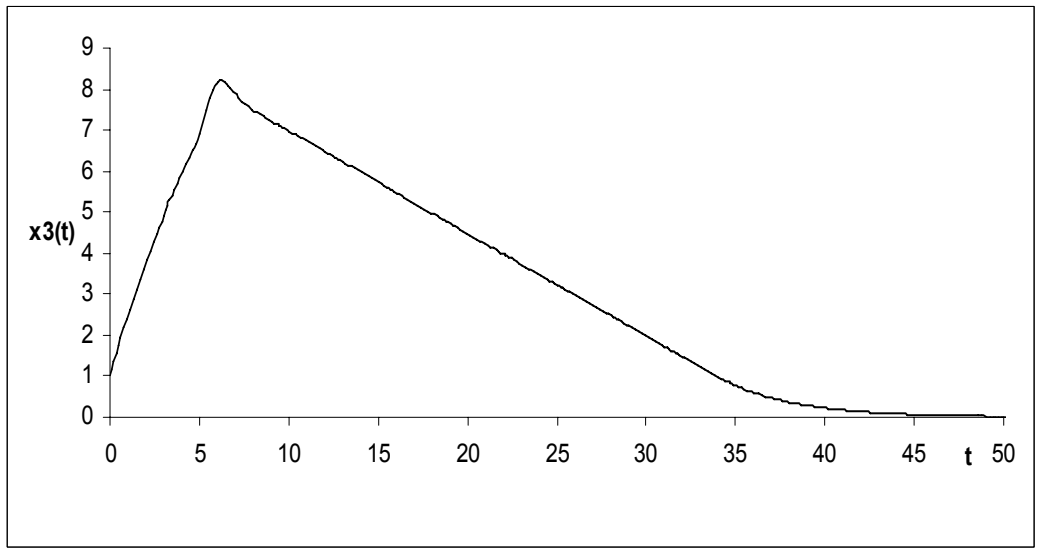

**Figure 1:** Time evolution of the state variables $x_1(t)$, $x_2(t)$ and $x_3(t)$ for the closed-loop system (4.1) with (3.18), where $k$ is defined by (4.6) with $r = 0.2$, $w(t) = \ln\left(\dfrac{2}{1+|\sin(t)|}\right)$ and initial condition
$$x_1(0) = x_2(0) = x_3(0) = 1$$

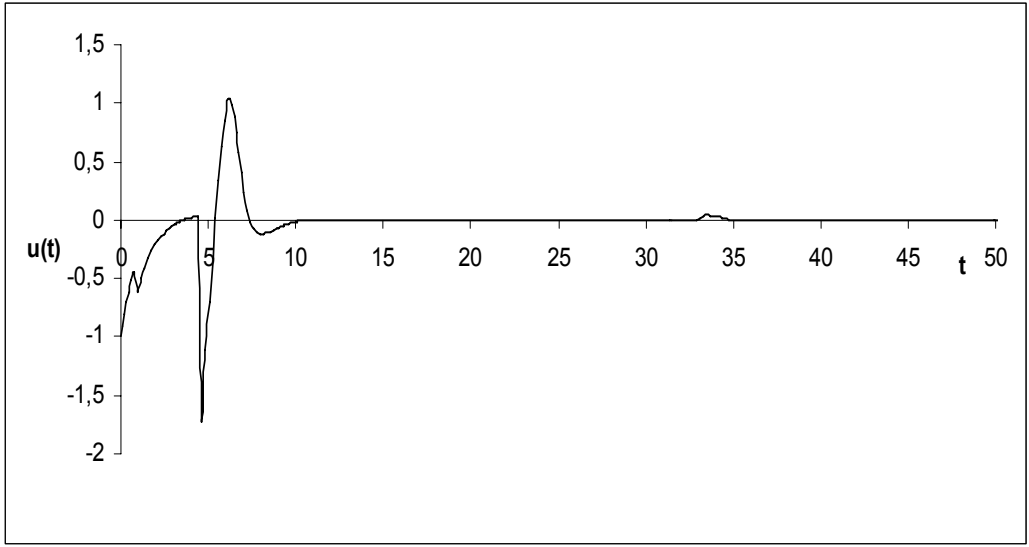

**Figure 2:** Time evolution of the input $u(t)$ for the closed-loop system (4.1) with (3.18), where $k$ is defined by (4.6) with $r = 0.2$, $w(t) = \ln\left(\dfrac{2}{1+|\sin(t)|}\right)$ and initial condition $x_1(0) = x_2(0) = x_3(0) = 1$

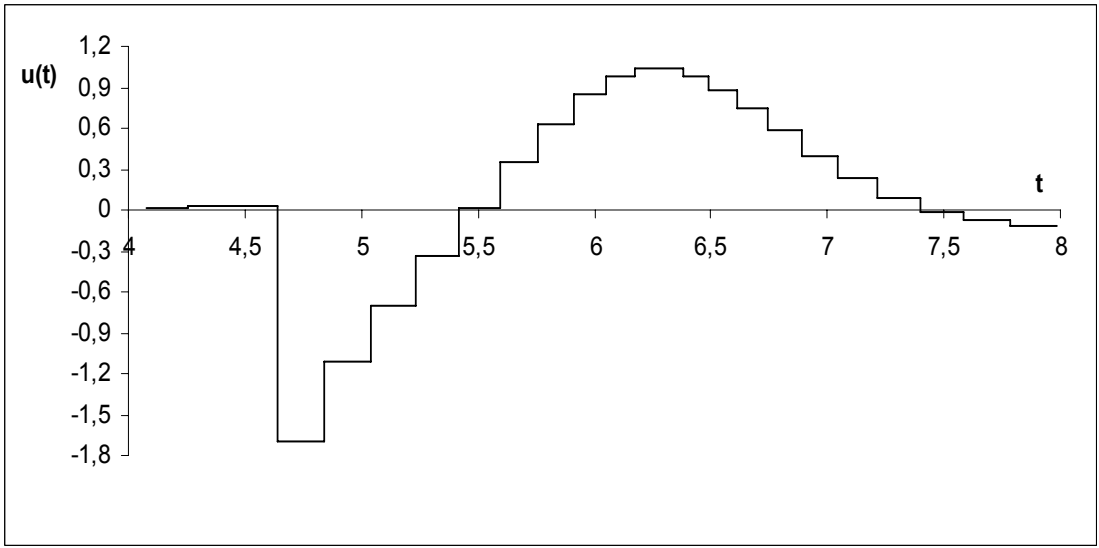

**Figure 3:** Time evolution of the input $u(t)$, $t \in [4,8]$ for the closed-loop system (4.1) with (3.18), where $k$ is defined by (4.6) with $r = 0.2$, $w(t) = \ln\left(\dfrac{2}{1+|\sin(t)|}\right)$ and initial condition $x_1(0) = x_2(0) = x_3(0) = 1$

We also assume the existence of a symmetric positive definite matrix $P \in \Re^{n \times n}$, a constant vector $p \in \Re^n$ and a constant $q > 0$ such that the matrix $P(A+bp') + (A'+pb')P$ is negative definite and such that

$$x'P(A+bp')x + x'Pf(d,x) \leq -q|x|^2, \text{ for all } (d,x) \in D \times \Re^n \tag{4.18}$$

Finally, we assume that

$$L_1 < \frac{qa_1|c'b|}{(1+|c|)a_2|Pb|} \quad (4.19)$$

where $c = -(A' + pb')^{-1}(0,\ldots,0,1)'$ and $a_2 \geq a_1 > 0$ are constants satisfying

$$a_1^2 x'Px \leq |x|^2 \leq a_2^2 x'Px, \text{ for all } x \in \Re^n \quad (4.20)$$

We show next that for every $\omega, R > 0$ with $a_2 R \leq \rho$ there exist constants $K, \tilde{r} > 0$ such that $0 \in \Re \times \Re^n$ is RGAS for the closed-loop system (4.15) with

$$\begin{aligned} u(t) &= \tilde{k}(x(\tau_i), y(\tau_i)), \, t \in [\tau_i, \tau_{i+1}) \\ \tau_{i+1} &= \tau_i + \tilde{r} \exp(-w(\tau_i)), \, \tau_0 = 0 \\ d(t) &\in D, w(t) \in \Re^+ \end{aligned} \quad (4.21)$$

where $\tilde{k} : \Re^n \times \Re \to \Re$ is defined by

$$\tilde{k}(x,y) := \begin{cases} p'x & \text{if } x'Px \geq R^2 \\ p'x - Kc'b \operatorname{sat}(\omega(y + c'x)) & \text{if } x'Px < R^2 \end{cases} \quad (4.22)$$

Indeed, this can be shown by a direct application of the Sampled-Data Forwarding Lemma with $k(x) := p'x$, arbitrary $\omega, R > 0$ with $a_2 R \leq \rho$, constants $M, K > 0$ satisfying

$$\frac{(1+|c|)L_1 a_2 R}{|c'b|^2} < K < \frac{qa_1 R}{|Pb||c'b|} \text{ and } M = \frac{K|c'b|\omega|Pb|}{(1+|c|)L_1} \quad (4.23)$$

and sufficiently small constant $\delta > 0$. Notice that by virtue of (4.18) and (4.20), inequality (3.3) is satisfied provided that $K|Pb||c'b| < qa_1 R$ (a direct consequence of (4.23)). Moreover, since $a_2 R \leq \rho$, it follows from (4.16) and (4.20) that inequality (3.4) holds provided that $(1+|c|)L_1 a_2 R < K(c'b)^2$ (a direct consequence of (4.23)). Finally, using the fact $a_2 R \leq \rho$ in conjunction with (4.16), (4.18), (4.20), we conclude that inequality (3.5) with $M$ as defined in (4.23) holds for sufficiently small $\delta > 0$ provided $(1+|c|)L_1|Pb| < q|c'b|$ (a direct consequence of (4.19)).

The only thing that remains to be shown is that hypothesis (H) of Lemma 3.1 holds with $k(x) := p'x$ and sufficiently small $r > 0$. By virtue of (4.18), we notice that for every $(x_0, d, w) \in \Re^n \times L_{loc}^\infty(\Re^+; D) \times L_{loc}^\infty(\Re^+; \Re^+)$, the solution of $\dot{x} = Ax + bu + f(d,x)$ with $x(0) = x_0$, $u(t) = p'x(\tau_i)$ satisfies the following differential inequality for almost all $t \in [\tau_i, \tau_{i+1})$:

$$\dot{V}(t) \leq -2q|x(t)|^2 + 2|x(t)||Pb||p||x(\tau_i) - x(t)| \quad (4.24)$$

where $V(t) = x'(t)Px(t)$, $\tau_{i+1} = \tau_i + r\exp(-w(\tau_i))$. By virtue of (4.17), it follows that the hypotheses of Lemma 2.3 hold for the absolutely continuous mapping $x : [\tau_i, \tau_{i+1}] \to \Re^n$. Using Lemma 2.3 and inequality (4.24), we conclude that for sufficiently small $r > 0$ there exists $\tilde{q} > 0$ such that the differential inequality $\dot{V}(t) \leq -\tilde{q}|x(t)|^2$ holds for almost all $t \geq 0$. Therefore, hypothesis (H) of Lemma 3.1 holds as well.

Notice that the case (4.15) includes systems, which are not necessarily feedforward systems. For example, the three-dimensional system:

$$\dot{x}_1 = k_1 d_1 x_1 + u$$
$$\dot{x}_2 = k_2 d_2 x_2 + x_1$$
$$\dot{y} = x_2 + d_3 g(x)$$
$$x = (x_1, x_2)' \in \Re^2, y \in \Re, d = (d_1, d_2, d_3)' \in [-1,1]^3, u \in \Re$$

where $k_1, k_2 > 0$ and $g \in C^1(\Re^2; \Re)$ with $g(0) = 0$, is not a feedforward system of the form (1.1). Inequalities (4.16), (4.17) hold for every $\rho > 0$ with $L_2 = L_1 = \max\{k_1, k_2, \max_{x \in B(0,\rho)} |\nabla g(x)|\}$. Moreover, inequality (4.18) holds with $P = \begin{bmatrix} 1 & 1+k_2 \\ 1+k_2 & (1+k_2)^2 + 1 \end{bmatrix}$, $p' = -[1 + S + k_2 \quad 1 + S(1+k_2)]$,

$S = \frac{1}{2} + k_1 + \frac{1}{2}(1+k_2)^2(k_2+k_1)^2$, $q = \dfrac{\sqrt{(1+k_2)^2 + 4} - 1 - k_2}{2 + 2k_2 + 2\sqrt{(1+k_2)^2 + 4}}$. ◁

## 5. Concluding Remarks

To construct a globally asymptotically stabilizing sampled-data feedback for feedforward systems subject to perturbations in the sampling schedule, we have developed the recursive sampled-data feedback synthesis tool—the Sampled-Data Forwarding Lemma. Assuming that a system (whose state is denoted by $x$) is stabilizable by sampled-data feedback, the Sampled-Data Forwarding Lemma guarantees the existence of sampled-data feedback stabilizer when the system is augmented by an additional state ($y$) in a cascade/feedforward configuration. Outside of the set $x'Px < R^2$ in the $(x,y)$ state space (where $P$ is an appropriate positive definite matrix and $R > 0$ is an appropriate constant) the feedback law uses the stabilizer for the $x$-subsystem in order to make $|x|$ small, whereas inside of the set $x'Px < R^2$ the feedback law uses a suitably saturated linear feedback law whose task is to drive $|y + c'x|$ below a prescribed small value $\omega^{-1}$ (where $c$ is an appropriate vector), while keeping $|x|$ small. Once both $|x|$ and $|y + c'x|$ are rendered small, the linear feedback law prevails in achieving exponential regulation to the origin. Robustness to perturbations in the sampling schedule is proved by quantifying the error between the current value of the state $x$ and its most recent sampled value, and by showing the smallness of this error for the closed-loop solutions provided all sampling periods fall uniformly below a sufficiently small "maximum allowable sampling period" (MASP).

Example 4.1 has shown that formulae (3.19), (3.20), (3.21) can be used in a straightforward way in order to design a globally stabilizing sampled-data feedback for an uncertain feedforward system of the form (1.1) under hypothesis (A2). However, the selection of the parameters $K_i > 0$, $R_i > 0$ ($i = 1, ..., n-1$) involved in formulae (3.19), (3.20), (3.21) is crucial for performance: low values for $K_i > 0$ will result in slow convergence of some state variables and high overshoot.

Example 4.1 has also demonstrated that the result of the present paper, in combination with the approach we introduced in [21], allows us to compensate any amount of actuation or sensing delay when controlling systems within the feedforward class using sampled-data controllers introduced in the present paper.

# Appendix

**Proof of Lemma 3.2:** Define:

$$\delta := -\max\{x'P(Ax + f(d,x,u) + bu) : (d,x) \in D \times \mathfrak{R}^n, x'Px = R^2, |u - p'x| \le K|c'b|\} > 0 \quad (A1)$$

The fact that $\delta$ as defined by (A1) is positive is a consequence of (3.3). Clearly, by virtue of continuity of the solution $x(t)$, there exists $\tau \in (0, \tau_1]$ such that $\max_{t \in [0,\tau]} x'(t)Px(t) < R^2$. The structure of system (3.7) guarantees that the solution of (3.7) exists for all times $\tau \in (0, \tau_1]$ with $\max_{t \in [0,\tau]} x'(t)Px(t) \le R^2$.

We prove by contradiction that $\max_{t \in [0,\tau_1]} x'(t)Px(t) < R^2$. We therefore assume that there exists $\tau \in (0, \tau_1]$ with $x'(\tau)Px(\tau) \ge R^2$. We define

$$T := \inf\{t \in [0,\tau_1] : x'(t)Px(t) \ge R^2\} \quad (A2)$$

and notice that $T \in (0, \tau_1]$. Definition (A2) and continuity of the solution $x(t)$ imply that $\max_{t \in [0,T]} x'(t)Px(t) = x'(T)Px(T) = R^2$. Define $V(t) = x'(t)Px(t)$ and $z(t) = y(t) + c'x(t)$. Notice that inequalities (3.2), (3.4) and the fact that $u(t) := p'x(0) - Kc'b\,\text{sat}(\omega z(0))$ for all $t \in [0,T)$ imply that the following inequality holds for almost all $t \in [0,T)$:

$$\begin{aligned}\dot{V}(t) &\le -2x'(t)P(Ax(t) + b(p'x(t) + v) + f(d(t), x(t), p'x(t) + v)) \\ &\quad + 2x'(t)P(f(d(t), x(t), p'x(0) + v) - f(d(t), x(t), p'x(t) + v)) + 2x'(t)Pbp'(x(0) - x(t))\end{aligned} \quad (A3)$$

where $v := -Kc'b\,\text{sat}(\omega z(0))$. The differential equation (3.7), in conjunction with $\max_{t \in [0,T]} x'(t)Px(t) \le R^2$ and $u(t) := p'x(0) - Kc'b\,\text{sat}(\omega z(0))$ imply that there exists a constant $S > 0$ such that $|\dot{x}(t)| \le S$ for almost all $t \in [0,T)$. Consequently, the following inequality holds for all $t \in [0,T]$:

$$|x(t) - x(0)| \le S t \quad (A4)$$

Using the facts that $f(d,x,u)$ is locally Lipschitz, $t \le T \le \tau_1 \le \tilde{r}$, $|v| \le K|c'b|$ and $\max_{t \in [0,T]} x'(t)Px(t) = x'(T)Px(T) = R^2$ in conjunction with definition (A1) and inequalities (A3), (A4), we guarantee the existence of a constant $L > 0$ such that the following inequality holds for almost all $t \in [0,T)$ sufficiently close to $T$:

$$\dot{V}(t) \le -\delta + L\tilde{r} \quad (A5)$$

If $\tilde{r} > 0$ is sufficiently small then inequality (A5) implies that $\dot{V}(t) \le -\delta/2$ for almost all $t \in [0,T)$ sufficiently close to $T$. This contradicts the assumption $\max_{t \in [0,T]} x'(t)Px(t) = x'(T)Px(T) = R^2$. The proof is complete. ◁

**Proof of Lemma 3.4:** Using the fact $c = -(A' + pb')^{-1}(0,...0,1)'$, we obtain $c'b \neq 0$ and $x_n + c'Ax = -c'bp'x$, for all $x \in \Re^n$. The previous equality in conjunction with the fact that $u(t) := p'x(0) - Kc'b\,\text{sat}(\omega z(0))$ for all $t \in [0, \tau_1)$ imply that the following differential equation holds for almost all $t \in [0, \tau_1)$:

$$\dot{z}(t) = c'bp'(x(0) - x(t)) + g(d(t), x(t), u(t)) + c'f(d(t), x(t), u(t)) - K(c'b)^2 \text{sat}(\omega z(0)) \tag{A6}$$

Define:

$$J := \max\{|g(d,x,u) + c'f(d,x,u)| : x'Px \leq R^2, d \in D, |u - p'x| \leq K|c'b|\} < K(c'b)^2 \tag{A7}$$

Using definition (A7), the fact that the mappings $f, g$ are locally Lipschitz and since $x'(t)Px(t) < R^2$ for all $t \geq 0$ (a consequence of Lemma 3.3), we obtain for all $t \in [0, \tau_1)$:

$$\begin{aligned}
&|g(d(t), x(t), u(t)) + c'f(d(t), x(t), u(t))| \leq \\
&|g(d(t), x(t), p'x(t) + v) + c'f(d(t), x(t), k'x(t) + v)| \\
&+ |g(d(t), x(t), p'x(0) + v) - g(d(t), x(t), p'x(t) + v)| \\
&+ |c||f(d(t), x(t), p'x(0) + v) - f(d(t), x(t), p'x(t) + v)| \\
&\leq J + L|x(t) - x(0)|
\end{aligned}$$

where $v = -Kc'b\,\text{sat}(\omega z(0))$ and $L > 0$ is an appropriate constant. Assuming that $K(c'b)^2 \omega \tilde{r} < 1$, integrating (A6) and exploiting the above inequality, we conclude that the following inequality holds for all $t \in [0, \tau_1]$:

$$|z(t)| \leq \left(1 - t\frac{K(c'b)^2 \omega}{\max(1, \omega|z(0)|)}\right)|z(0)| + t\left[Q \max_{0 \leq s \leq t}|x(s) - x(0)| + J\right] \tag{A8}$$

where $Q > 0$ is an appropriate constant. The differential equations (3.7) in conjunction with $\max_{t \in [0,\tau_1]} x'(t)Px(t) \leq R^2$ and $u(t) := p'x(0) - Kc'b\,\text{sat}(\omega z(0))$ imply that there exists a constant $S > 0$ such that $|\dot{x}(t)| \leq S$ for almost all $t \in [0, \tau_1)$. Consequently, inequality (A4) holds for all $t \in [0, \tau_1]$. Combining (A4), (A8), we can conclude that the following inequality holds for all $t \in [0, \tau_1]$:

$$|z(t)| \leq \left(1 - t\frac{K(c'b)^2 \omega}{\max(1, \omega|z(0)|)}\right)|z(0)| + QSt^2 + Jt \tag{A9}$$

The above inequality in conjunction with inequality (3.4) (which implies that $J < K(c'b)^2$), shows that the following implications hold for sufficiently small $\tilde{r} > 0$:

- If $\omega|z(0)| \geq 1$ then $|z(t)| \leq |z(0)|$ for all $t \in [0, \tau_1]$
- If $\omega|z(0)| \leq 1$ then $|z(t)| \leq \omega^{-1}$ for all $t \in [0, \tau_1]$

It follows that $|z(t)| \leq \max(|z(0)|, \omega^{-1})$ for all $t \in [0, \tau_1]$. Using induction, it can be shown that:

$$|z(t)| \leq \max(|z(\tau_i)|, \omega^{-1}) \text{ for all } t \geq \tau_i, \; i = 0,1,2,... \tag{A10}$$

Moreover, inequality (A9) shows that the following implication holds for sufficiently small $\tilde{r} > 0$:

$$\text{If } \omega|z(\tau_i)| \geq 1 \text{ then } |z(\tau_{i+1})| \leq |z(\tau_i)| - G(\tau_{i+1} - \tau_i) \tag{A11}$$

where $G > 0$ is an appropriate constant.

Implication (A11) shows that (3.9) holds with $T(z) := \frac{\max(0, \omega|z| - 1)}{\omega G} + \tilde{r}$. Indeed, we prove this by contradiction. Suppose that there exists $t > \frac{\max(0, \omega|z(0)| - 1)}{\omega G} + \tilde{r}$ with $|z(t)| > \omega^{-1}$. Let $m \in Z^+$ be the largest integer with $\tau_m \leq t < \tau_{m+1}$. By virtue of (A10) we conclude that $|z(\tau_m)| > \omega^{-1}$. Moreover, since $\tau_{m+1} \leq \tau_m + \tilde{r}$ it follows that $\tau_m > \frac{\max(0, \omega|z(0)| - 1)}{\omega G}$. Estimate (A10) shows that $|z(\tau_i)| > \omega^{-1}$ for all $i = 0, \ldots, m$. Implication (A11) gives $|z(\tau_m)| \leq |z(0)| - \tau_m G$, which is a contradiction.

The proof is complete.   ◁

**Proof of Lemma 3.5:** By virtue of Lemma 3.3 and Lemma 3.4 the solution of (3.7) satisfies $x'(t)Px(t) < R^2$, $|z(t)| \leq \omega^{-1}$ for all $t \geq 0$. Let $t \geq 0$ be a time where $V(t) = \frac{M}{2}z^2(t) + \frac{1}{2}x'(t)Px(t)$ is differentiable. Let $m \in Z^+$ be the largest integer with $\tau_m \leq t < \tau_{m+1}$. Using (A6) we obtain:

$$\dot{V}(t) = S_1(t) + S_2(t) \tag{A12}$$

where

$$\begin{aligned}S_1(t) &:= -MK(c'b)^2 \omega z^2(t) - Kc'b\omega z(t)x'(t)Pb \\&+ x'(t)P(Ax(t) + bp'x(t) + f(d(t), x(t), p'x(t) - Kc'b\omega z(t))) \\&+ Mz(t)g(d(t), x(t), p'x(t) - Kc'b\omega z(t)) \\&+ Mz(t)c'f(d(t), x(t), p'x(t) - Kc'b\omega z(t))\end{aligned} \tag{A13}$$

$$\begin{aligned}S_2(t) &:= [Mz(t)c'b + x'(t)Pb]p'(x(\tau_m) - x(t)) \\&+ Mz(t)[g(d(t), x(t), p'x(\tau_m) - Kc'b\omega z(\tau_m)) - g(d(t), x(t), p'x(t) - Kc'b\omega z(t))] \\&+ (x'(t)P + Mz(t)c')[f(d(t), x(t), p'x(\tau_m) - Kc'b\omega z(\tau_m)) - f(d(t), x(t), p'x(t) - Kc'b\omega z(t))] \\&+ Kc'b\omega[Mc'b z(t) + x'(t)Pb](z(t) - z(\tau_m))\end{aligned} \tag{A14}$$

Notice that inequality (3.5) implies that

$$S_1(t) \leq -\delta z^2(t) - \delta |x(t)|^2 \tag{A15}$$

Moreover, since the mappings $f, g$ are locally Lipschitz and $x'(t)Px(t) < R^2$, $|z(t)| \leq \omega^{-1}$ for all $t \geq 0$, it follows that the hypotheses of Lemma 2.3 hold on the interval $[\tau_m, \tau_{m+1}]$ for the absolutely continuous map $(z(t), x(t))$ for appropriate constants $Q, G$. Therefore, for sufficiently small $\tilde{r} > 0$, there exists $\Gamma > 0$ such that the following inequality holds for all $t \in [\tau_m, \tau_{m+1}]$

$$\max\left(|x(t)-x(\tau_m)|,|z(t)-z(\tau_m)|\right) \le \Gamma \tilde{r}\,|(x(t),z(t))| \tag{A16}$$

Using the facts that the mappings $f,g$ are locally Lipschitz and $x'(t)Px(t) < R^2$, $x'(\tau_m)Px(\tau_m) < R^2$, $|z(\tau_m)| \le \omega^{-1}$, $|z(t)| \le \omega^{-1}$ in conjunction with inequality (A16) and definition (A14), we obtain:

$$S_2(t) \le \Gamma \tilde{r} q |z(t)|^2 + \Gamma \tilde{r} q |x(t)|^2 \tag{A17}$$

for certain appropriate constant $q > 0$. Selecting $\tilde{r} > 0$ sufficiently small and using inequalities (A15), (A17), we can conclude that (3.10) holds with $\mu = \frac{\delta}{2}$. The proof is complete. ◁

**Proof of Lemma 3.6:** Since $P \in \Re^{n \times n}$ is a symmetric positive definite matrix, there exist constants $a_2 \ge a_1 > 0$ satisfying

$$a_1^2 x'Px \le |x|^2 \le a_2^2 x'Px, \text{ for all } x \in \Re^n \tag{A18}$$

Since $P(A+bp') + (A'+pb')P$ is negative definite there exists a constant $q > 0$ such that

$$x'P(A+bp')x \le -q|x|^2, \text{ for all } x \in \Re^n \tag{A19}$$

Let $C < \dfrac{qa_1}{|Pb||c'b|}$ with $C \in (0,1]$ (otherwise arbitrary) and $R^* > 0$ such that

$$Q(R^*)a_2 R^* < \dfrac{q|c'b|}{(1+|c|)\left(\dfrac{\lambda q}{1+|p|}+|Pb|\right)\left(1+|p|+a_2^{-1}C|c'b|\right)+(1+|p|)|P||c'b|} \tag{A20}$$

$$Q(R^*)a_2 R^* < \dfrac{C(c'b)^2}{(1+|c|)\left((1+|p|)a_2+C|c'b|\right)} \tag{A21}$$

$$Q(R^*)a_2 R^* < \dfrac{qa_1 - C|Pb||c'b|}{(1+|p|)|P|a_2+|P|C|c'b|} \tag{A22}$$

where $Q(R) := L\left((1+|p|)a_2 R + R|c'b|\right)$. We claim that there exists a constant $\delta > 0$ such that (3.3), (3.4), (3.5) hold with arbitrary $\omega > 0$, $R \in (0, R^*)$, $K = CR$, $M = \dfrac{C|c'b|\omega}{(1+|c|)Q(R)}\dfrac{|P|Q(R)a_2 R + |Pb|}{(1+|p|)a_2 + C|c'b|}$ for the case $Q(R) > 0$ and $M = \dfrac{CR|Pb|^2 \omega}{4q} + 1$ for the case $Q(R) = 0$. Indeed, using (3.17), (A18), (A19), in conjunction with the fact that $C \le 1$, we conclude that conditions (3.3), (3.4) with $K = CR$ are satisfied provided that

$$(1+|p|)|P|Q(R)a_2^2 R + (|Pb|+|P|Q(R)a_2 R)C|c'b| < qa_1 \tag{A23}$$

$$(1+|c|)(1+|p|)Q(R)a_2^2 R + (1+|c|)Q(R)a_2 CR|c'b| < C(c'b)^2 \tag{A24}$$

Inequalities (A23), (A24) are direct consequences of (A21), (A22) and the fact that $R \leq R^*$. Finally, after tedious calculations and using (3.17), (A18), (A19) in conjunction with the fact that $C \leq 1$, we conclude that condition (3.5) is satisfied provided that there exists $\delta > 0$ such that the following inequality holds for all $(x, z) \in \mathfrak{R}^n \times \mathfrak{R}$:

$$\left(M(1+|c|)(1+|p|)Q(R)a_2 R + K|P|c'b|\omega Q(R)a_2 R + K|c'b\|Pb|\omega + M(1+|c|)Q(R)K|c'b|\right)|x\|z|$$
$$- MK(c'b)^2 \omega z^2 - \left(q - (1+|p|)|P|Q(R)a_2 R\right)|x|^2 \leq -\delta|z|^2 - \delta|x|^2$$

The existence of sufficiently small $\delta > 0$ such that the above inequality holds is a direct consequence of (A20), the facts that $K = CR$, $R \leq R^*$ and the selection $M = \dfrac{C|c'b|\omega}{(1+|c|)Q(R)} \dfrac{|P|Q(R)a_2 R + |Pb|}{(1+|p|)a_2 + C|c'b|}$ for the case $Q(R) > 0$ and $M = \dfrac{CR|Pb|^2 \omega}{4q} + 1$ for the case $Q(R) = 0$. The proof is complete. ◁